\documentclass{article}
\usepackage{twocolextabst}

\usepackage{jlcode}
\usepackage{hyphenat}
\usepackage{times}
\usepackage{soul}
\usepackage{url}
\usepackage[hidelinks]{hyperref}
\usepackage[utf8]{inputenc}
\usepackage[small]{caption}
\usepackage{graphicx}
\usepackage{amsmath}
\usepackage{booktabs}
\usepackage{subcaption}
\usepackage{tikz}
\usepackage{amssymb}
\usetikzlibrary{calc}
\usetikzlibrary{shapes.geometric}
\usetikzlibrary{decorations.markings}

\usepackage{xargs}
\usepackage{todonotes}
\newcommandx{\unsure}[2][1=]{\todo[inline,linecolor=red,backgroundcolor=red!25,bordercolor=red,#1]{#2}}
\newcommandx{\change}[2][1=]{\todo[inline,linecolor=blue,backgroundcolor=blue!25,bordercolor=blue,#1]{#2}}
\newcommandx{\info}[2][1=]{\todo[inline,linecolor=green,backgroundcolor=green!25,bordercolor=green,#1]{#2}}
\newcommandx{\improvement}[2][1=]{\todo[inline,linecolor=purple,backgroundcolor=purple!25,bordercolor=purple,#1]{#2}}
\newcommandx{\thiswillnotshow}[2][1=]{\todo[disable,#1]{#2}}

\title{Compositional Scientific Computing with Catlab and SemanticModels}

\author{
Micah Halter$^1$\footnote{Corresponding Author}\and
Evan Patterson$^2$\and
Andrew Baas$^1$\And
James P. Fairbanks$^1$
\affiliations
$^1$Georgia Tech Research Institute, Atlanta, GA USA\\
$^2$Stanford University, Stanford, CA USA
\emails
micah.halter@gtri.gatech.edu
}

\newcommand{\semanticModels}{\nohyphens{\texttt{SemanticModels.jl}}}
\newcommand{\catlab}{\nohyphens{\texttt{Catlab.jl}}}
\newcommand{\petri}{\nohyphens{\texttt{Petri.jl}}}
\newcommand{\diffeq}{\nohyphens{\texttt{DifferentialEquations.jl}}}

\begin{document}

\maketitle

\begin{abstract}
Scientific computing is currently performed by writing domain specific modeling frameworks for solving special classes of mathematical problems. Since applied category theory provides abstract reasoning machinery for describing and analyzing diverse areas of math, it is a natural platform for building generic and reusable software components for scientific computing. We present \catlab{}, which provides the category-theoretic infrastructure for this project, together with \semanticModels{}, which leverages this infrastructure for particular modeling tasks. This approach enhances and automates scientific computing workflows by applying recent advances in mathematical modeling of interconnected systems as cospan algebras. 
\end{abstract}

\thiswillnotshow{We should probably write this from the perspective of our new Brand Algebraic Julia?\\
- I think we could keep it pre-rebrand (split SM and Catlab), and use the TAC full Catlab paper as the announcement for that and also includes the SemanticModels cospan algebra/solving as features in Catlab and the different bridges to other modeling frameworks as separate julia packages.}

\thiswillnotshow{sounds good for this abstract. Maybe for the talk we present the new brand. Evan and I decided to announce the new brand at juliacon, which will be before this talk is presented. So let's use algebraic julia ecosystem as the combined brand and describe the two components seperately.}

\paragraph{Introduction}
The AlgebraicJulia ecosystem is a suite of software libraries for building scientific computing software in Julia based on categorical logic. The core components are \catlab{}, which provides the data structures and algorithms for symbolic computing in a categorical setting, and \semanticModels{}, which facilitates several metamodeling tasks by representing models
as open systems using decorated cospans~\cite{1605.08100,1906.05443}. This paradigm allows metamodeling tasks
to be rigorously defined in the language of category theory and then applied across diverse scientific domains.
\semanticModels{} leverages these categorical definitions to both to facilitate the
metamodeling process and to generate solvers for the new models.

\paragraph{GATs and Partially Symbolic Computing}
\catlab{} implements a computer algebra system for the fragment of dependent type theory known as \emph{generalized algebraic theories }(GATs) \cite{cartmell1986}. These theories provide a formal language for defining algebraic structures that includes all the familiar algebraic objects such as groups, rings, and modules, as well as order structures such as preorders, posets, and lattices, and  categorical structures such as categories, monoidal categories, bicategories and double categories, and quantum computing's ZX calculus. In the case of monoidal categories, including those with extra algebraic gadgets like comonoids or internal homs, \catlab{} provides three syntactical systems for representing morphisms. The first is a point-free formula notation like $(f \otimes g) \cdot h$ which is stored as an abstract syntax tree (AST) like the LISP program \nohyphens{\texttt{(compose (otimes f g) h)}}. The second is a wiring diagram or string diagram notation illustrated in Figure~\ref{fig:string_diagram}.
The third is a program syntax, shown in Figure~\ref{fig:program_syntax}, that provides a familiar interface for programmers to express morphisms.
The goal of \catlab{} is to bridge the divide between symbolic computing, where programs reason about mathematical models, and scientific computing, where programs compute numerical simulations of mathematical models.

\begin{figure}[htb!]
\centering
\begin{subfigure}[b]{0.9\linewidth}
\centering
\begin{tikzpicture}[unit length/.code={{\newdimen\tikzunit}\setlength{\tikzunit}{#1}},unit length=4mm,x=\tikzunit,y=\tikzunit,semithick,box/.style={rectangle,draw,solid,rounded corners},outer box/.style={draw=none},wire/.style={draw,postaction={decorate},decoration={markings,mark=at position 0.5 with {\node[anchor=south] {#1};}}}]
  \node[outer box,minimum width=10\tikzunit,minimum height=7\tikzunit] (root) at (0,0) {};
  \node[box,minimum size=2\tikzunit] (n3) at (-2,1.5) {$f$};
  \node[box,minimum size=2\tikzunit] (n4) at (-2,-1.5) {$g$};
  \node[box,minimum width=2\tikzunit,minimum height=5\tikzunit] (n5) at (2,0) {$h$};
  \path[wire=$A$] ($(root.west)+(0,1.5)$) to[out=0,in=-180] (n3.west);
  \path[wire=$B$] ($(root.west)+(0,-1.5)$) to[out=0,in=-180] (n4.west);
  \path[wire=$B$] (n3.east) to[out=0,in=-180] ($(n5.west)+(0,1.5)$);
  \path[wire=$C$] (n4.east) to[out=0,in=-180] ($(n5.west)+(0,-1.5)$);
  \path[wire=$D$] (n5.east) to[out=0,in=180] (root.east);
\end{tikzpicture}
\caption{Catlab can draw string diagrams representations of morphisms using TikZ, Graphviz, or \texttt{Compose.jl}.}%
\label{fig:string_diagram}
\end{subfigure}
\begin{subfigure}[b]{.9\linewidth}
\begin{lstlisting}[language=Julia]
p = @program MonoidalCategory (a::A, b::B) begin
    return h(f(a), g(b))
end
\end{lstlisting}
\caption{\catlab{} supports a traditional program syntax for defining morphisms in a syntax familiar to programmers.}\label{fig:program}
\label{fig:program_syntax}
\end{subfigure}
\caption{The expression ${(f \otimes g) \cdot h}$ depicted in the different syntactic representations supported by \catlab{}.}%
\label{fig:catlab_syntax}
\end{figure}

\paragraph{Model Composition}

We can construct a scientific open system as a model, ${N\colon X \to Y}$, where $X$ and $Y$ are sets of ``input'' and ``output'' variables respectively. We can then represent this system 
as a cospan ${f\colon X \to N, g\colon Y \to N}$ in some category $\mathcal{C}$ drawn ${X \to N \gets Y}$~\cite{1502.00872}.
In order to represent open scientific systems, we define this category
$\mathcal{C}$ as a symmetric monoidal bicategory with cospans as morphisms, following Courser's example ~\cite{1605.08100}. Therefore, $\mathcal{C}$ not only has composition,
defined as ${G \circ F\colon X \to Z}$ for cospans
${F\colon X \to Y}$ and ${G\colon Y \to Z}$, but also a monoidal product, defined as the bifunctor
${\otimes\colon \mathcal{C}\times\mathcal{C} \to \mathcal{C}}$ where if
${F\colon X_1 \to Y_1}$ and ${G\colon X_2 \to Y_2}$ then ${F \otimes
G\colon X_1 \otimes X_2 \to Y_1 \otimes Y_2}$. In this category, composition
represents the serial combination of two open systems. This
composition is equivalent to the first system executing given an
input and passing its output to the second system, which then executes and
produces the final output. The monoidal product
of two open systems represents the models running independently in parallel,
where the resulting system's input is the union of the original models' inputs,
the output is the union of the original models' outputs, and the model
simulates the two systems executing simultaneously. When $\mathcal{C}$ is $\mathbf{FinSet}$, this monoidal product is the coproduct or disjoint union, meaning we can represent objects as tuples of variables with $X \otimes Y$ formed by concatenation.

\begin{figure}[htb!]
\includegraphics[width=\linewidth]{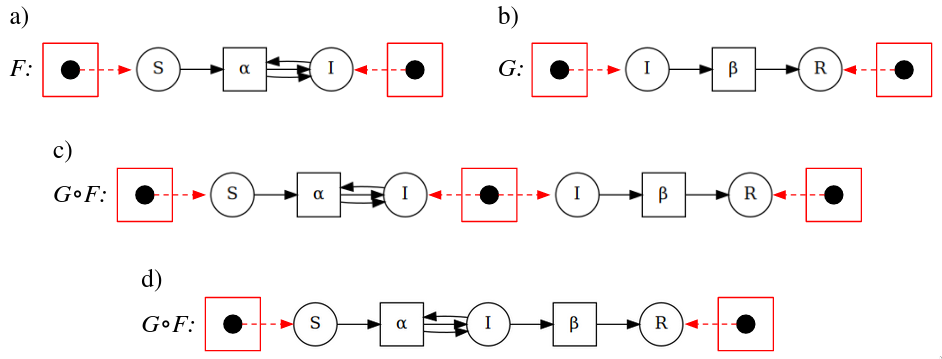}
\caption{Open Petri models (a)~${F\colon X \to M \gets Y}$, (b)~${G\colon
Y \to N \gets Z}$, and (c,d)~${G \circ F\colon X \to M+_{Y}N \gets Z}$.}%
\label{fig:composition}
\end{figure}

Using these two operations as an interface, we can build a framework where
scientists can easily create more complex models by manipulating and composing basic
models. Consider the open Petri
net models\footnote{Petri net: a simple graphical model that consists of states,
denoted as circles which can store tokens or particles, and transitions that move those particles between states, denoted as
squares connected to states by arrows. Petri nets are frequently used in a variety of fields such as systems
biology and epidemiology~\cite{Gorrieri_2017,9783540412175}.}\textsuperscript{,}\footnote{Recent work by Baez and Master
defines a symmetric monoidal bicategory of Petri Nets~\cite{1808.05415}.} in
Figure~\ref{fig:composition}.  Model $F$ (Figure~\ref{fig:composition}a) has
states $S$ and $I$ which represent susceptible and infected populations respectively, and
$\alpha$ which represents an interaction between a susceptible and an infected person
that results in two infected people.  Similarly Model $G$
(Figure~\ref{fig:composition}b) has states $I$ and $R$ which represent infected
and recovered populations respectively, and transition $\beta$ which represents an infected person recovering.
Figures~\ref{fig:composition}c and~\ref{fig:composition}d depict the model
$G \circ F$, an SIR epidemiology model~\cite{Bahi_Jaber_2003}.

\begin{figure}[htb!]
\centering
\includegraphics[width=.7\linewidth]{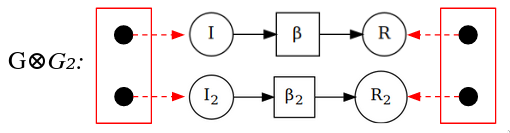}
\caption{The monoidal product of two instances of the model $G$ results in ${Y \otimes Y \to N \otimes N \gets Z \otimes Z}$.}%
\label{fig:tensor}
\end{figure}

Figure~\ref{fig:tensor} demonstrates the monoidal product of two instances of
model $G$ (Figure~\ref{fig:composition}b). Model $G \otimes G$ can be
interpreted as two non-interacting infected populations recovering from their
respective illnesses separately.

\paragraph{Model Comparison}

This method of building models compositionally extends beyond just model construction to support more complex metamodeling tasks, which can be viewed as higher order operations on models. For example, the categorical representation of 
models provides an intuitive method of model comparison.
Scientists make decisions about the governing laws of the phenomena of interest with significant uncertainty. In this case they need to rapidly compare the results of simulating many similar models. By describing models with morphism formulas, we can rapidly compare models that are formed by substituting different components into the same expression.

\begin{figure}[htb!]
\centering
\includegraphics[width=\linewidth]{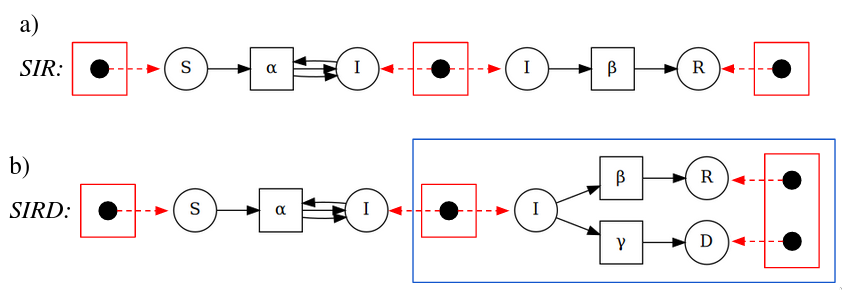}
\caption{Comparison of the open system representations of (a) an SIR epidemiological model expressed as the composition $G\circ F$ and (b) an SIRD epidemiological model expressed as the composite $H\circ F$. In both models the relationship between between susceptible and infectious individuals is the same.}%
\label{fig:comparison}
\end{figure}

Figure~\ref{fig:comparison} demonstrates this method of model comparison by presenting two models. Figure~\ref{fig:comparison}a presents an open system representation of an SIR epidemiology model which can be defined as $G \circ F$ where $F$ is the infection model from Figure~\ref{fig:composition}a and $G$ is the recovery model from Figure~\ref{fig:composition}b. Figure~\ref{fig:comparison}b presents an SIRD epidemiology model defined as $H \circ F$ where $H$ is an alternative recovery model, boxed in blue, with a new transition, $\gamma$, that represents a spontaneous reaction where a person moves from being infected ($I$) to dead ($D$). With this approach of comparing model substructures, we not only identify the different states and transitions between the models, but we recognize that they are instances of the infection model $F$ composed with different recovery models. In this way, the hierarchical formula representation gives us a way to describe the relationships between models.

\thiswillnotshow{(This has been addressed, but I have left as reference to think about) In the algebraic approach we generally think of objects as the important thing and we want the morphisms as ways to map between objects preserving the algebraic structure that defines our category. In our approach the objects are spaces where things can happen and the morphisms are processes on those spaces. That is pretty different. We want to harmonize these approaches and talk about the similarity between two morphisms based on \emph{having the same structure with different substructures.}}

\paragraph{Implementation}

\semanticModels{} \cite{semanticmodelsjl} implements these model composition and comparison ideas to create a generalizable and extensible scientific metamodeling framework.
While standard approaches to software development with category theory (namely functional programming) treat programs as morphisms in a $\mathbf{Set}$-like category, we implement category theory as a library within the Julia programming language~\cite{doi:10.1137/141000671}. In this paradigm, the scientific models are represented as morphisms in a category and implemented as data structures within a program, in contrast to viewing the program itself as a morphism in a category. 

\semanticModels{} depends on \catlab{} \cite{catlabjl}, an applied category theory framework in Julia, to automate these metamodeling tasks for scientists, such as comparing and composing categorically defined models.
Additionally, \semanticModels{} utilizes features in \catlab{} to interface these categorically defined models with existing modeling frameworks, such as \petri{}~\cite{petrijl} and \diffeq{}~\cite{diffeqjl}, to calculate simulations and solutions.
This allows categorical construction of
models with an intuitive method of formally describing new models, and provides
a method of generating solvers for these new models with minimal software development effort.

\begin{figure}[htb!]
\centering
\includegraphics[width=.7\linewidth]{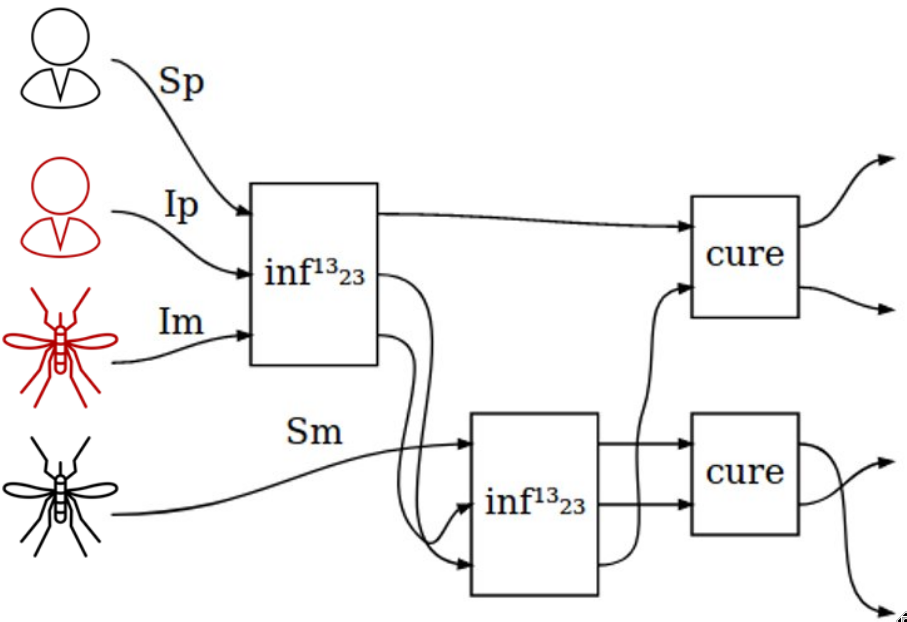}
\caption{Wiring diagram representation of the composition of epidemiology infection and recovery models to create a new model that simulates a complex dual infection system. $S_p$ and $I_p$ represent ``susceptible'' and ``infected'' people, respectively, while $S_m$ and $I_m$ represent ``susceptible'' and ``infected'' mosquitoes.}%
\label{fig:malaria}
\end{figure}

\begin{figure}[htb!]
\centering
\includegraphics[width=\linewidth]{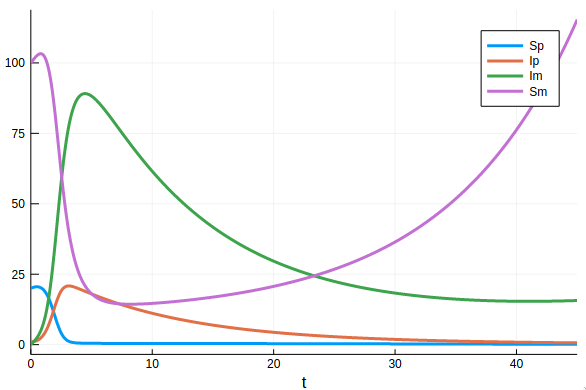}
\caption{A solution to the system defined in Figure~\ref{fig:malaria}, automatically generated using \semanticModels{}. Parameters of the model are chosen by the scientist by choosing parameters for each box in the string diagram, these correspond to choosing reaction rates in the corresponding Petri Net. In this way, a parameterized scientific model corresponds to a model of a theory in the language of categorical logic.}%
\label{fig:result_graph}
\end{figure}

For example, Figure~\ref{fig:malaria} shows a wiring diagram produced by \semanticModels{} where each box is a decorated cospan and the input and output wires show the applied composition operations. Using simple building blocks that simulate infection and recovery, we can easily build a complex system that models the dynamics between two populations that are affected by malaria spread in a dual infectious manner. \semanticModels{} can then interpret that wiring diagram and automatically generate a simulation that takes as input parameters and initial conditions and produces a solution trajectory as shown in Figure~\ref{fig:result_graph}.

\paragraph{Conclusion}
\catlab{} and the associated ecosystem of tools support
scientific progress by enabling the rapid adaptation and extension of models from existing works to new scientific phenomena and engineering problems. By implementing and building upon the recent developments in decorated cospans, we are able to create a general-purpose metamodeling framework. This framework can use known models and domain-specific concepts as building blocks for modeling novel systems, can easily compare models both structurally and numerically, and can automatically generate simulations without time consuming and error prone programming effort.

\bibliographystyle{vancouver}
\bibliography{refs}
\end{document}